\documentclass[10pt]{article}
\usepackage[utf8]{inputenc}
\usepackage{mathtools,
        mathrsfs,amssymb,amsthm, amsmath,
        bm}
\usepackage{hyperref}
\usepackage{url}
\usepackage[a4paper, left = 1in, right = 1in, top =1in, bottom = 1in]{geometry} 
\usepackage{ dsfont }   
\usepackage[british]{babel}

\theoremstyle{definition}
\newtheorem*{rem}{Remark}

\def\[#1\]{\begin{align*}#1\end{align*}}

\newcommand{\R}{\mathbb{R}}
\newcommand{\ceq}{\coloneqq}
\newcommand{\eps}{\varepsilon}
\newcommand{\df}{\mathop{}\!\mathrm{d}}

\begin{document}

\title{A Remark on Lebesgue Criterion}
\author{Yu-Lin Chou\thanks{Yu-Lin Chou, Institute of Statistics, National Tsing Hua University, Hsinchu 30013, Taiwan,  R.O.C.; Email: \protect\url{y.l.chou@gapp.nthu.edu.tw}.}}
\date{}
\maketitle

\begin{abstract}
We remark a variant of the existence part of the fundamental theorem
of calculus, which, together with the Lebesgue differentiation theorem,
constitute a new proof that every Riemann-integrable function on a
compact interval having limit everywhere on the interior is almost
everywhere continuous with respect to Lebesgue measure. The proof
is intended as a new connection between Lebesgue differentiation theorem
and Lebesgue criterion of Riemann integrability.\\

{\noindent \textbf{Keywords:}} almost everywhere continuity; fundamental theorem of calculus; Lebesgue criterion;  Lebesgue differentiation theorem; Riemann integrability

{\noindent \textbf{MSC 2020:}} 26A42; 28A15
\end{abstract}

\section{Introduction}
A characterization for Riemann integrability is Lebesgue criterion,
i.e. the assertion that a bounded function on a compact interval is
Riemann-integrable if and only if the function is continuous almost
everywhere modulo Lebesgue measure. As well-known and well-documented,
the more involved part of Lebesgue criterion (for and only for a direct
proof) would be the only-if part, which, even with measure theory
employed here and there, usually takes some nontrivial preliminary
observations. The only-if part of Lebesgue criterion is thus more
“interesting”. 

We obtain a new proof of a weakened form of the only-if part\footnote{We recall that boundedness follows from Riemann integrability.}
of Lebesgue criterion: \textit{If $[a, b] \subset \R$, if $f: [a,b] \to \R$ is Riemann-integrable, and if $f$ has limit everywhere on $]a,b[$, then $f$ is continuous almost everywhere with respect to Lebesgue measure.} 
The main observation is a natural variant of the existence part of
the fundamental theorem of calculus — the assertion that every continuous
function on a compact interval has a primitive, the indefinite integral
of the function being a suitable choice. The variant seems to be easily
unnoticed, but, as will be shown, takes only a simple proof. Moreover,
our proof of the proposition connects Lebesgue criterion with Lebesgue
differentiation theorem in a new direction. 

\section{Proof}
Let $f$ satisfy the assumptions of our proposition. Then the Riemann
integrability of $f$ implies that $\int_{a}^{x}f$ exists for every
$x \in [a,b]$. Since $f$ has limit everywhere on $]a,b[$, let $l_{x} \ceq \lim_{y \to x}f(y)$
for every $x \in ]a,b[$. For convenience, write $\partial_{x}$ for
the usual differential operator. 

We claim that $\partial_{x}\int_{a}^{x}f = l_{x}$ for every $x \in ]a,b[$.
If $x \in ]a,b[$, and if $h > 0$, then \[
\bigg| h^{-1}\bigg( \int_{a}^{x+h}f - \int_{a}^{x}f \bigg) - l_{x} \bigg|
&= \bigg| h^{-1}\int_{x}^{x+h}f - l_{x} \bigg|\\
&= \bigg| h^{-1}hl_{x} - l_{x} + h^{-1}\int_{x}^{x+h}(f(t) - l_{x}) \df t \bigg|\\
&= \bigg| h^{-1}\int_{x}^{x+h} ( f(t) - l_{x}) \df t \bigg|\\
&\leq h^{-1}\int_{x}^{x+h}|f(t) - l_{x}| \df t.
\]Given any $\eps > 0$, there is some $\delta > 0$ such that $0 < h < \delta$
implies that $|f(t) - l_{x}| < \eps$ for all $t \in ]x, x+h]$; so
the last term above is $< \eps$ for all $0 < h < \delta$. The case
where $h < 0$ follows for the same reason after repeating the analysis
for the expression $|h|^{-1}\int_{x+h}^{x}f$. This proves the claim.

On the other hand, we have $\lim_{h \to 0+}(2h)^{-1}\int_{x - h}^{x+h}f = \partial_{x}\int_{a}^{x}f$
for every $x \in ]a,b[$. Since $f$ is Riemann-integrable, the extension
of $f$ to $\R$ by assigning $0$ to every point of $\R \setminus [a,b]$
is Lebesgue-integrable; the Lebesgue differentiation theorem (e.g.
Theorem 7.2, Wheeden and Zygmund \cite{wz}) asserts in particular
that $(2h)^{-1}\int_{x-h}^{x+h}f \to f(x)$ as $h \to 0+$ for almost
every $x \in ]a,b[$ (modulo Lebesgue measure). But then $l_{x} = \lim_{y \to x}f(y) = f(x)$
for Lebesgue-almost all $x \in ]a,b[$, and the desired almost everywhere
continuity follows. \qed

\begin{rem}
It is not \textit{a priori} obvious whether a Riemann-integrable
function on a compact interval should even have one-sided limit everywhere.
Call a function on $[a,b]$ a step function if and only if there is
some partition $\{ x_{0}, \dots, x_{n} \}$ of $[a,b]$, including
the endpoints $a,b$, such that the function is constant on each $]x_{i}, x_{i+1}[$.
Since it can be shown (e.g. Theorem 7.6.1, Dieudonn{\'e} \cite{d})
that a function on a compact interval has one-sided limit everywhere
if and only if the function is the uniform limit of some sequence
of step functions, and since not every Riemann-integrable function
on a compact interval is the uniform limit of some sequence of step
functions, being Riemann-integrable alone does not even ensure having
one-sided limit everywhere. \qed 
\end{rem}


\begin{thebibliography}{1} 

\bibitem{d}
Dieudonn{\'e}, J. (1969). 
\textit{Foundations of Modern Analysis}.
Academic Press.

\bibitem{wz}
Wheeden, R.L. and Zygmund, A. (2015).
\textit{Measure and Integral: An Introduction to Real Analysis},
second edition.
CRC Press.

\end{thebibliography}
\end{document}